\documentclass[12pt]{article}
\usepackage{amssymb}
\begin{document}
\title{Weak covering and the 
tree property}
\author{Ralf-Dieter Schindler\footnote{The author
gratefully acknowledges a grant of the Deutsche 
Forschungsgemeinschaft (DFG)}}

\maketitle

{\abstract Suppose that there is no transitive model of $ZFC \ +$
there is a strong cardinal, and let $K$ denote the core model. It is
shown that if $\delta$ has the tree property then $\delta^{+K} =
\delta^+$ and $\delta$ is weakly compact in $K$.} 

\bigskip
Let $\lnot L^{strong}$ denote the assumption that there is no 
transitive (set or proper class sized)
model of $ZFC \ +$ there is 
a strong cardinal. We have shown in \cite{WC} (cf.
Theorem 5) that 
$K^c$ correctly computes successors of
weakly compact cardinals, provided that $\lnot L^{stong}$ holds.
Here, $K^c$ is the countably complete core model below a strong
cardinal. In fact, a straightforward
adaptation of the argument given in \cite{WC}
also yields that $K$, the true core model below a strong
cardinal, has the same weak covering property. (Schimmerling had
earlier proved this for the core model below a Woodin cardinal,
cf. \cite{schsteel}.
See \cite{jensen}, 
\cite{diss}, or \cite{zeman} on the theory of $K^c$ and $K$.)

The purpose of this note is to prove a more general result, namely the
following 

\bigskip
{\sc Theorem.} {\it If $\lnot L^{strong}$ holds 
then $\delta^{+K} = \delta^+$
for every cardinal $\delta$ having the tree property,
and such a $\delta$ is weakly compact in $K$.}

\bigskip
A cardinal $\delta$ is said to have the tree property if there is no
Aronszajn $\delta$-tree, i.e., if every tree of height $\delta$ whose
levels all have size less than $\delta$ has a cofinal branch. 
It is elementary that $\omega$ has the tree property but
$\omega_1$ does not. By \cite{mitchell}, $\omega_2$ (and in fact
for example every double successor) may have the tree
property, given the consistency of a weakly compact cardinal.
 
It is more difficult to show that $\omega_2$ and $\omega_3$ may 
simultaneously have the tree property (cf. \cite{conseq}). In fact,
Foreman and Magidor have shown in unpublished work that at least a
Woodin cardinal seems to be neccessary to do the job.
Our Theorem may be viewed as saying that one reason for the difficulty
in forcing $\omega_2$ and $\omega_3$ to 
simultaneously have the tree property is that
weak covering not only holds for singular and weakly compact cardinals
(cf. \cite{jensen}, \cite{WC}, viz. 
\cite{mitsch}, \cite{schsteel}), but
also for cardinals admitting the tree property. 

We expect our Theorem
to generalize to much higher core models, but we do not at the
moment know how to prove it for Steel's core model below one Woodin
cardinal, say.
 
The proof of the Theorem 
consists in applying the following lemmata to
the core model. The first of them is a simple observation building 
upon a classical insight of Jensen. For its statement, 
recall Jensen's principle
$\square^\star_\kappa$, cf. \cite{fine}, p. 283.
 
\bigskip 
{\sc Lemma 1.1.} {\it Let $\delta$ have the tree property. 
Suppose $W$ to
be an inner model such that for some $W$-cardinal $\kappa$ we have 
$W \models 2^\kappa = \kappa^+ \wedge \square^\star_\kappa$. Then
$\delta \not= \kappa^{+W}$. In particular, $\delta$ is inaccessible
in any inner model $W$ in which $GCH$ and $\square^\star_\kappa$ hold
for all $\kappa < \delta$.}

\bigskip
{\it Proof.} By \cite{fine}
p. 283, inside $W$, using $2^\kappa = \kappa^+$ and
$\square^\star_\kappa$ we can construct a special Aronszajn
$\kappa^+$-tree. So if $\delta = \kappa^{+W}$ then in $V$ there is
an Aronszajn $\delta$-tree. Contradiction!

\hfill $\square$ (Lemma 1.1)

\bigskip
To state the second, and main, 
lemma, let us introduce the following terminology.
Let ${\cal H}$, ${\tilde {\cal H}}$ be two models of the same type. We
call an elementary embedding $\sigma \colon {\cal H} \rightarrow
{\tilde {\cal H}}$ $\kappa$-complete (for a cardinal $\kappa \geq
\aleph_1$) iff for every elementary 
$\tau \colon {\bar {\cal H}} \rightarrow {\tilde {\cal H}}$ with
$Card({\bar {\cal H}}) < \kappa$ there is an elementary 
$\pi \colon {\bar {\cal H}} \rightarrow {\cal H}$ such that
$\sigma \circ \pi (x) = \tau (x)$ for all $x \in {\bar {\cal H}}$ 
such that $\tau(x) \in ran(\sigma)$. If $\sigma$ is $\aleph_1$-complete
then it is also called countably complete.

In this situation, in particular, if ${\cal H} = (M;\in,...)$ is a
transitive structure then ${\tilde {\cal H}} = ({\tilde M};E,...)$
is well-founded and we thus may and will
identify ${\tilde {\cal H}}$ with its 
transitive collapse. Note that if we had
an infinite decreasing $E$-sequence $... \ E \ x_1 \ E \ x_0$ then we
could choose $\tau \colon {\bar {\cal H}} \rightarrow {\tilde {\cal
H}}$ with $\{ x_n \colon n < \omega \} \subset ran(\tau)$ and the
countable completeness of $\sigma$ would give some $\pi \colon
{\bar {\cal H}} \rightarrow {\cal H}$ and $... \in \pi \circ 
\tau^{-1}(x_1) \in \pi \circ \tau^{-1}(x_0)$. Contradiction!      

\bigskip
{\sc Lemma 1.2.} {\it Let $\delta$ have the tree property.
Let $W$ be an inner model such that $\delta$ is (strongly)
inaccessible in $W$, $H =
(H_{\delta^{+}})^W$ ($=$ the set of all sets
in $W$ being hereditarily $\leq \delta$ in $W$)
has cardinality $\delta$, and 
$cf(\delta^{+W}) = \delta$. Then there is a $\delta$-complete
$\sigma \colon H \rightarrow {\tilde H}$ such that
$sup \ \sigma$''$\delta^{+W} < On \cap {\tilde H}$.}

\bigskip
{\it Proof.} To commence, we note that every $X \subset \delta^{+W}$
of cardinality $< \delta$ can be covered by some $Y \in W$ of 
cardinality $< \delta$. [Let w.l.o.g. $sup(X) > \delta$. As 
$cf(\delta^{+W}) = \delta$, there is $g \in W$, $g \colon
\delta \rightarrow sup(X)$ bijective. But $\delta$ is
regular, as it has the tree property, so $\theta =
sup(g^{-1}$''$X) < \delta$, and
$Y = g$''$\theta \in W$ is such that $X \subset Y$.] 
 
Now let $F \colon \delta \rightarrow H$ be bijective. By the
previous paragraph and $\delta$'s being inaccessible in $W$
we may pick
$(A_\xi : \xi < \delta)$ such that for all ${\bar \xi} <
\xi < \delta$ we have $\delta \in A_\xi \in W$, $A_{\bar \xi} \subset
A_\xi$, $Card(A_{\bar \xi}) < 
Card(A_\xi) < \delta$, and $F$''$\xi \subset A_\xi$.
For every $n < \omega$ let 
$h_n \colon \omega \times H \rightarrow H$, $h_n \in K$, 
be some $\Sigma_{n+1}$ Skolem
function for $H$ being definable over $H$. 
For every $X \in {\cal P}(H) \cap K$ let us write $h[X]$ for 
$\bigcup_{n < \omega} \ h_n$''$(\omega
\times X)$ where $X \subset H$, noting that $h[X] \prec H$. 
Trivially,
$H = \bigcup_{\xi < \delta} \ h[A_\xi]$. 

We now let $T$ consist of all $(\xi , \eta)$ such that
$\xi < \delta <
sup(\delta^{+W} \cap h[A_\xi]) 
< \eta < \delta^{+W}$. Note that by $cf(\delta^{+W}) =
\delta$ for every $\xi < \delta$ there are $\delta$ many
$\eta$'s with $(\xi , \eta) \in T$. We consider $T$ as being ordered
by setting $(\xi , \eta) \leq_T ({\tilde \xi} , {\tilde \eta})$
iff $\xi \leq {\tilde \xi}$ and there is 
$\pi \colon h[A_\xi \cup \{ \eta \}] \rightarrow
h[A_{\tilde \xi} \cup \{ {\tilde \eta} \}]$ induced by
the requirements $\pi \upharpoonright A_\xi = id$ and
$\pi(\eta) = {\tilde \eta}$.

Set $[\xi , \eta] = \{ (\xi , {\tilde \eta}) \in T \colon
(\xi , \eta) \leq_T (\xi , {\tilde \eta}) \leq_T (\xi , \eta) \}$,
and let $T^\star$ be the set of all such $[\xi , \eta]$'s.
Obviously, $\leq_T$ induces a tree ordering $\leq_{T^\star}$
on $T^\star$. In fact, $( T^\star , \leq_{T^\star})$ can
be checked to be a $\delta$-tree. [The $\xi$'s level of
$T^\star$ consists of nodes of the form $[\xi , \eta]$ for
some $\eta$. Now suppose that this level had cardinality
$\delta$, say $\{ [\xi , \eta^i] \colon i < \delta \}$ with $[\xi ,
\eta^i] \not= [\xi , \eta^j]$ for $i < j < \delta$ are its
nodes. Using $(2^{Card(A_\xi)+ \aleph_0})^W < \delta$ 
and the pigeonhole principle
we may then find $i < j < \delta$ such that $(\xi , \eta^i)
\leq_T (\xi , \eta^j) \leq_T (\xi , \eta^i)$. Contradiction!]

Now let $b$ be any cofinal branch thru $T^\star$ given by the 
tree property of $\delta$. Let us write $\pi_{\xi,{\tilde \xi}}$ for 
$\pi \colon h[A_\xi \cup \{ \eta \}] \rightarrow
h[A_{\tilde \xi} \cup \{ {\tilde \eta} \}]$ given by $[\xi , \eta]
\leq_{T^\star} [{\tilde \xi} , {\tilde \eta}] \in b$. Let 
$({\tilde H} , \pi_{[\xi , \eta],b})$ be the direct limit of the
system $(h[A_\xi \cup \{ \eta \}] , \pi_{\xi,{\tilde \xi}})$.
We may define $\sigma \colon H \rightarrow {\tilde H}$ by sending
$x \in H$ to that thread having eventually constant value $x$. 

It is now easy to check that $\sigma$ is as desired.
Let $\tau \colon {\bar H} \rightarrow {\tilde H}$ be elementary such
that ${\bar H}$ has cardinality $< \delta$. 
Using the regularity of $\delta$,
$ran(\tau) \subset ran(\pi_{[\xi,\eta],b})$ for some $[\xi,\eta]
\in b$. Then $\pi \colon {\bar H} \rightarrow H$ is well-defined and
elementary where we set $\pi(x) = h_n(m,\gamma)$ for $\tau(x) =
\pi_{[\xi,\eta],b}(h_n(m,\gamma))$, $n, m < \omega$, 
$\gamma \in A_\xi \cup \{ \eta \}$.
Moreover, $\tau(x) \in ran(\sigma)$ means that $\tau(x) =
\pi_{[\xi,\eta],b}(h_n(m,\gamma))$ for some $n, m < \omega$ and some 
$\gamma \in A_\xi$; but then $\sigma \circ \pi(x) = \tau(x)$.
 
We also have that the thread given by the
$\eta$'s for $[\xi,\eta] \in b$, $\xi < \delta$, is above every thread
having constant value $\zeta$ for any $\zeta < \delta^{+W}$ which
implies that $ sup \ \sigma$''$On \cap H =
sup \ \sigma$''$\delta^{+W} < On \cap {\tilde H}$. 

\bigskip
\hfill $\square$(Lemma 1.2)

\bigskip
Requiring also that $\delta \cap h[A_\xi]$ is an ordinal, we could have
arranged that $\sigma \upharpoonright \delta = id$. Moreover, by
replacing the requirement $\xi < \delta < sup(\delta^{+W} \cap
h[A_\xi]) < \eta < \delta^{+W}$ by $\xi < \delta$ and $sup(\delta \cap
h[A_\xi]) < \eta < \delta$ we can arrange that in fact $\delta$ is the
critical point of $\sigma$ (then $T^\star$ resembles the tree
constructed in the proof of Lemma 2 in
\cite{JSL}), and we obtain the following

\bigskip
{\sc Lemma 1.3.} {\it Let $\delta$ have the tree property. Let $H$ be
a transitive model of $ZFC \ - \{ Powerset \}$ such
that $\delta$ is (strongly) inaccessible in $H$, $H$ has
cardinality $\delta$, and $cf(On \cap H) = \delta$. 
Then there is a $\delta$-complete $\sigma
\colon H \rightarrow {\tilde H}$ with critical point $\delta$.}  

\bigskip
Before now turning toward the proof of the Theorem let us remark that
as a matter of fact if $\lnot L^{strong}$ holds then $cf(\kappa^{+K})
\geq Card(\kappa)$ for every $\kappa \geq \aleph_2$. Jensen has shown
this in \cite{jensen} (Theorem 7) for $\kappa \geq \aleph_3$, but the
proof of \cite{mitsch} in fact yields this slight strengthening. 

\bigskip
{\it Proof} of the Theorem. Let us fix a cardinal $\delta$ having the
tree property. Let us assume that $\lnot L^{strong}$ holds (in
particular $K$, the true core model below one strong cardinal exists),
however $\delta^{+K} < \delta^+$. We shall derive a contradiction.

Jensen has shown that $\square_\kappa$ holds in $K$ for every
$K$-cardinal $\kappa$ (cf. \cite{zeman}, where it is even shown that
$K \models \square$, i.e., global square holds in $K$). 
In particular, $\square^\star_\kappa$ holds
everywhere and so by Lemma 1.1 $\delta$ is inaccessible in $K$.
By the above remark, $cf(\delta^{+K}) \geq
\delta$, and so $cf(\delta^{+K}) = \delta$. 
Hence by Lemma 1.2, setting $H =
(H_{\delta^+})^K$, the set of all sets 
in $K$ being hereditarily $\leq \delta$
in $K$, we may choose a countably complete
$$\sigma \colon H \rightarrow {\tilde H}.$$

\bigskip
{\sc Claim.} ${\tilde H}$ is a mouse.

\bigskip
{\it Proof.} This is a standard argument.
Suppose not, and let ${\cal I}$ be a putative 
iteration of
${\tilde H}$ with a last ill-founded model. Let $\theta$ be large
enough and let $\tau^\star \colon {\bar V} \rightarrow V_\theta$ be
elementary such that ${\bar V}$ is countable and transitive, and
$\{ {\tilde H} , {\cal I} \} \subset 
ran(\tau^\star)$. By absoluteness, 
${\bar {\cal I}} = 
\tau^{\star -1}({\cal I})$ is a (countable) putative iteration of
${\bar H} = \tau^{\star -1}({\tilde H})$ with a last ill-founded
model. On the other hand, setting $\tau = \tau^\star
\upharpoonright {\bar H}$, by the countable completeness of $\sigma$
there is an elementary $\pi \colon {\bar H} \rightarrow H$ with
$\sigma \circ \pi(x) = \tau(x)$ for all $x \in {\bar H}$ such that 
$\tau(x) \in ran(\sigma)$, so that ${\bar {\cal I}}$ can be copied to
give an iteration ${\cal I}_0$ of $H$. But $H$ is a mouse, so the last
model of ${\cal I}_0$ is well-founded. But then the last model of
${\bar {\cal I}}$ is well-founded, too. Contradiction!

\bigskip
\hfill $\square$ (Claim)

\bigskip
We now let $W = Ult(K;\sigma)$, the ultrapower of $K$ using $\sigma$
as an extender (cf. for example \cite{mitschsteel} \S 2.5), 
and we let $${\tilde \sigma} \colon
K \rightarrow W$$ denote the associated ultrapower map. 
We have that $\sigma(\delta)^{+W} = 
{\tilde \sigma}(\delta)^{+W} = sup \ {\tilde
\sigma}$''$\delta^{+K} = sup \ \sigma$''$\delta^{+K}$. 
Moreover, $W$ is a universal weasel. [$K$ correctly computes cofinally
many in $On$ successors, and so does $W$. But this implies the
universality of $W$ by $\lnot L^{strong}$.] 

We may now coiterate $W$ with ${\tilde H}$. As $W$ is universal, 
there can be no truncation on the ${\tilde H}$-side of the
coiteration. As $\lnot L^{strong}$, $\delta$ is not overlapped in $K$,
i.e., there is no $K$-measurable $\mu < \delta$ such that $\mu$
is a strong cardinal in ${\cal J}^K_\delta$. [Otherwise we have found
a transitive model of $ZFC$ + there is a strong cardinal.]
So $\sigma(\delta)$ is not overlapped neither in $W$ nor in
${\tilde H}$. Hence the coiteration is above $\sigma(\delta)$ on the
${\tilde H}$-side in the sense that the critical points of all
extenders used are $\geq \sigma(\delta)$. 
But $\sigma(\delta)$ is the largest cardinal of
${\tilde H}$, so that in fact ${\tilde H}$ is not moved at all in this
coiteration.

Hence we know that ${\tilde H}$ is an initial segment of an
iterate $W^\star$ of $W$, where the iteration giving $W^\star$ 
is above $\sigma(\delta)$, too. 
But now, by the property of
$\sigma$ given by Lemma 1.2, 
$$\sigma(\delta)^{+W} = sup \ \sigma^{\prime \prime}\delta^{+K} 
< On \cap {\tilde H} \leq \sigma(\delta)^{+W^\star} \leq
\sigma(\delta)^{+W}.$$ Contradiction!

We have thus shown that $\delta^{+K} = \delta^+$, and we 
are left with having to verify that $\delta$ is weakly compact in
$K$. For this it suffices to show that in $K$, for
all $S \subset {\cal P}(\delta)$ with $Card(S) \leq \delta$ there is a
uniform $\delta$-complete filter deciding $S$.

Let $S \subset {\cal J}^K_\alpha$ with $\alpha 
< \delta^{+K} = \delta^+$, and pick
$\beta > \alpha$,
$\beta < \delta^+$, such that ${\cal J}^K_\beta \models Card(\alpha)
= \delta$, $cf(\beta) = \delta$, and 
${\cal J}^K_\beta \prec {\cal J}^K_{\delta^+}$. 
Pick $f \in {\cal J}^K_\beta$ with $f \colon \delta
\rightarrow {\cal P}(\delta) \cap {\cal J}^K_\alpha$ bijective. Using
Lemma 1.3 we get a $\delta$-complete $$\sigma \colon
{\cal J}^K_\beta \rightarrow {\tilde K}$$ 
with critical point $\delta$. Define $U = \{ X \in {\cal
P}(\delta) \cap {\cal J}^K_\alpha \colon \delta \in \sigma(X) \}$.

For $X \in {\cal P}(\delta)$ we have
$X \in U$ iff there is $\xi < \delta$ with $X = f(\xi)$ and
$\delta \in \sigma(f)(\xi)$. But $f$ may be coded by a subset 
of $\delta$ in ${\cal J}^K_\beta$, so $f \in {\tilde K}$. 

Thus $U \in {\tilde K}$, and $U$ is coded
by a subset $A_U$ of $\delta$ in ${\tilde K}$. Now the coiteration of
$K$ with ${\tilde K}$ is above $\delta$ on both sides by the same
argument that gave above that the coiteration of $W$ with ${\tilde
H}$ is above $\sigma(\delta)$ on both sides. Moreover, the
coiteration is simple
on the $K$-side by the universality of $K$. 
But this implies that $A_U \in K$, and so $U \in K$.
Finally, the $\delta$-completeness of $\sigma$ implies that $U$ is
$\delta$-complete in $K$, and clearly $U$ is uniform.

\bigskip
\hfill $\square$ (Theorem)

\bigskip
{\sc Mathematisches Institut

Uni Bonn

Beringstr. 4

53115 Bonn, Germany}

{\tt rds@rhein.iam.uni-bonn.de}

\bigskip
{\sc Math Dept.

UCB

Berkeley, CA 94720, USA}

{\tt rds@math.berkeley.edu}


\bigskip
\begin{thebibliography}{99}
\bibitem{conseq} Abraham, U., {\it Aronszajn trees on $\aleph_2$
and $\aleph_3$}, Ann. Pure Appl. Logic {\bf 24} (1983), pp. 213 - 230.
\bibitem{fine} Jensen, R. B., {\it The fine structure of the
constructible hierarchy}, Ann. Math. Logic {\bf 4} (1972), pp. 229 -
308.
\bibitem{jensen} ---------, {\it The core model for non-overlapping
extender sequences}, handwritten notes.
\bibitem{mitchell} Mitchell, W. J., {\it Aronszajn trees and the
independence of the transfer property}, Ann. Pure Appl. Logic {\bf 5}
(1972), pp. 21 - 46.
\bibitem{mitsch} ---------, and Schimmerling, E.,
{\it Covering without countable closure}, Math. Research Letters {\bf
2} (1995), pp. 595 - 609.
\bibitem{mitschsteel} ---------, ---------, and Steel, J., {\it
The covering lemma up to a Woodin cardinal}, Ann. Pure Appl. Logic.
\bibitem{schsteel} Schimmerling, E., and Steel, J., 
{\it The maximality of the core model}, preprint.
\bibitem{diss} Schindler, R.-D. {\it The core model up to one strong
cardinal}, Ph. D. Thesis, Bonn, 1996.
\bibitem{WC} ---------, {\it Weak covering at large cardinals}, 
Math. Logic Quarterly {\bf 43} (1997), pp. 22 - 28.
\bibitem{JSL} ---------, {\it Successive weakly compact or singular
cardinals}, Journal Symb. Logic, to appear.
\bibitem{zeman} Zeman, Martin, {\it The core model for non-overlapping
extender sequences and its applications}, Ph. D. Thesis, Berlin, 1997.
\end{thebibliography}
\end{document}